\documentclass{amsart}

\usepackage{graphicx}

\newtheorem{theorem}{Theorem}[section]
\newtheorem{lemma}[theorem]{Lemma}
\newtheorem{proposition}[theorem]{Proposition}

\theoremstyle{corollary}
\newtheorem{corollary}[theorem]{Corollary}

\theoremstyle{definition}

\theoremstyle{remark}

\numberwithin{equation}{section}

\begin{document}

\title[]{Strong asymptotic behavior of multi-orthogonal polynomials associated with a queueing model}

\author{U. Fidalgo}
\address{Department of Mathematics and Statistics,
Case Western Reserve University, Cleveland, Ohio 43403}
\email{ufx6@case.edu}


\subjclass[2000]{Primary 60H10, 60J35; Secondary 41X00}

\date{\today}


\keywords{Multiple orthogonal polynomials}

\begin{abstract} We describe the strong asymptotic behavior of  type I and type II multiple orthogonal polynomials which were used to give an integral expression for a transition probability function corresponding to a queueing models that has a bulk service admitting batches with a fixed size of $m\in \mathbb{N}$ customers.
\end{abstract}

\maketitle                                        

\section{Introduction}\label{introduction}

In  \cite{F} we present a connection between a bulk queueing model and sequences of muti-orthogonal polynomials. This work concerns with the strong asymptotic behavior of theses sequences of polynomials involeved. The mentioned queueing model (considered in  \cite{B, Med}) describes queues with a fixed number $m\in \mathbb{N}=\{1,2,3,\ldots\}$ of admitted costumers, while the individuals to be served arrive one by one. Let $X(t)$  be the random variable that counts the waiting customers. Assume that both customers and servers follow two Poisson processes with constant parameters $\lambda>0$ and  $\mu>0$, respectively. In terms of Kendall's notation \cite{Ke} (also used in \cite{Med}) this system is denoted by $M/M(m,m)/1$. It corresponds to a stationary Markov process whose transition probability functions has form
\[
P_{n,r}(t)=Pr \left\{X(t+s)=r | X(s)=i\right\}, \quad (n,r) \in \mathbb{Z}_+\times \mathbb{Z}_+=\mathbb{Z}_+^2,
\]
where $\mathbb{Z}_+=\mathbb{N}\cup \{0\}=\{0,1,2,\ldots\}$.  We gave (see  \cite{F}) an integral formula for the functions $P_{n,r}$, $(n,r)\in \mathbb{Z}_+^2$ in terms of the elements of the sequences of polynomials $\displaystyle \left\{Q_n\right\}_{n\in \mathbb{Z}_+}$ and sequences of vector polynomials $\displaystyle \left\{{\bf q}_r=(q_{0,r},\ldots,q_{m-1,r})\right\}_{r\in \mathbb{Z}_+}$
\[
P_{n,r}(t)=\int e^{x t} Q_n(x) \sum_{j=0}^{m-1} q_{j,r}(x) \, d \sigma_j(x),
\]
where  $\sigma_j,$ $j=0,1,\ldots, m$, represent certain distributions that we describe bellow. 

Given three real numbers $\alpha$, $\beta$, and $\theta$ with $\alpha <\beta$ and $\theta \in [0,\, 2 \pi]$, we denote by $\displaystyle \left[\alpha,\, \beta\right]\exp i \theta$ the segment in the complex plane that results after rotating the interval $[\alpha,\beta]$, an angle of $\theta$ with counterclockwise orientation, and center at the endpoint $\alpha$. This notation could also be extended to the case when $[\alpha, \infty)$. Let us consider the starlike set  
\begin{equation}\label{Sigma0}
\Sigma_0=\bigcup_{k=0}^m \left[-\lambda-\mu, -\lambda-\mu + \frac{m+1}{m}\left(\frac{\mu\lambda^m}{m}\right)^{1/(m+1)} \right] \exp \frac{2 \pi i k}{m+1}.
\end{equation}

Let $\delta_{\zeta}$ denote the Dirac delta distribution which is a probability measure with a point of mass at  $\zeta \in \mathbb{R}$. Following definitions in \cite[Chapter 6]{Rud}, for each $j\in \mathbb{Z}_+$, $\delta_{\zeta}^{(j)}$ denotes the $j$th distribution derivative of $\delta_{\zeta}$. This means that any  $j$th differentiable function $f$ at $\zeta$ whose derivatives $f^{(k)}$, $k=0,1,\ldots,j$, are integrable in $\mathbb{R}$, satisfies that
\begin{equation}\label{deltaderivative}
\int f(z) \, \delta_{\zeta}^{(j)}(z)\,  d z=(-1)^{j}f^{(j)}(\zeta).
\end{equation}
We can extend  (\ref{deltaderivative}) for any differentiable function $f$ in an open interval that contains  $\zeta$, replacing $f$ by a proper test function expressed by a piecewise definition which coincides with $f$ on such interval. Notice that when $j\in \mathbb{N}$ the distribution $\delta_{\zeta}^{(j)}$ is not a measure. The system $(\sigma_0, \sigma_1, \ldots, \sigma_m)$ contains $m$ distributions that have a type of  Radon Nikodym derivative
\[
\frac{d \sigma_j}{d x}=\rho_j (x) \rho(x)+\frac{(-1)^{j-1}}{\lambda^{j}(j-1)!}\left(\delta_{-\lambda}-\delta_{-\lambda-\mu}\right)^{(j-1)}, \quad j=0,1,\ldots, m-1,
\]
agreeing $\rho_0 \equiv 1$ and $\displaystyle \frac{\delta^{(-1)}}{(-1)!}=0$. The expression $\rho_j (x) \rho(x)$ is  a function defined on the starlike set where, fixed $k \in \{0,1,\ldots, m\}$ and each point 
\[
x_k=-\mu-\lambda+ t \exp \frac{2 \pi i k}{m+1} \in \Sigma_0, \quad 0\leq t \leq  \frac{m+1}{m}\left(\frac{\mu\lambda^m}{m}\right)^{\frac{1}{m+1}},
\]
the following symmetry holds
\[
\rho_j \left(x_k\right) \rho \left(x_k\right)= \exp \left(-\frac{2 \pi i j}{m+1} \right)\rho_j \left(t-\mu-\lambda\right) \rho \left(t-\mu-\lambda\right),
\]
being $\displaystyle \rho_j \left(t-\mu-\lambda\right) \rho \left(t-\mu-\lambda\right)$ a weight. The functions $\rho, \rho_1, \ldots, \rho_{m-1}$ are expressed explicitly in terms of a solution $\omega_0$ of the algebraic equation (see \cite{AKS,F})
\begin{equation}\label{omegaequationQ}
\lambda\, \omega^{m+1}-(z+\lambda+\mu)\, \omega^m+\mu=0,
\end{equation}
as follows
\[
\rho_j(x)=\frac{1}{2 \pi i} \left(\frac{1}{\omega_{0+}^j(x)}-\frac{1}{\omega_{0-}^j(x)}\right), \quad x\in \Sigma, \quad j=1,\ldots,m,
\]
where for each $\displaystyle t\in (0, a)$ with $a=\frac{m+1}{m} \left(\frac{\mu\lambda^m}{m}\right)^{1/(m+1)}$
\[
\lim_{h\to 0, h>0} \omega_0\left(-\mu-\lambda+(t\pm ih)\exp \left(\frac{2 \pi i k}{m+1}\right)\right)
\]
\[
=\omega_{0\pm}\left(-\mu-\lambda+t\exp \left(\frac{2 \pi i k}{m+1}\right)\right), \quad k=0,1, \ldots,m.
\]
The branch chosen is such that  $\omega_0$ is analytic in $\mathbb{C}\setminus \Sigma_0$ (we denote $\omega_0 \in \mathcal{H}\left(\mathbb{C}\setminus \Sigma_0\right)$) satisfying
\begin{equation}\label{asimptoticomega0}
\omega_0(z)= \lambda z+\lambda+\mu +\mathcal{O}\left(\frac{1}{z}\right) \quad \mbox{as} \quad z \to \infty.
\end{equation}
Observe that to first components of the system of distribution, $\sigma_0$ and $\sigma_1$, are measures. 

The sequences of $\displaystyle \left\{Q_n\right\}_{n\in \mathbb{Z}_+}$ and $\displaystyle \left\{{\bf q}_r=(q_{0,r},\ldots,q_{m-1,r})\right\}_{r\in \mathbb{Z}_+}$ are constructed by the following two higher order three terms recurrence relations
 \begin{equation}\label{recurrenceQ}
(\lambda+\mu+x)Q_n(x)=\lambda Q_{n+1}(x) +\mu Q_{n-m}(x), \quad n \geq m,  
\end{equation}
with initial conditions
\[
Q_n(x)=\frac{1}{\lambda^n}(\lambda+x)^{n},\quad n=0,\ldots, m,
\]
and
\begin{equation}\label{preccurencelatinq}
\begin{array}{l}
(\lambda+x) {\bf q}_r(x)=\mu{\bf q}_{r+m}(x)+\lambda {\bf q}_{r-1}(x), \quad  r \in \{0,\ldots, m-1\},\\ \\
(\lambda+\mu+x) {\bf q}_r(x)=\mu{\bf q}_{r+m}(x)+\lambda {\bf q}_{r-1}(x), \quad  r \in \{m,m+1,\ldots\}, \\ \\ \mbox{with initial coonditions} \\ \\
\displaystyle {\bf q}_{-1}=(0,\ldots,0), {\bf q}_{0}=\left(1,0,\ldots,0\right),  \\ \\
{\bf q}_{1}=\left(0,1,0,\ldots,0\right), \ldots, {\bf q}_{m-1}=\left(0,\ldots,, 0,1\right).
\end{array}
\end{equation}

In  \cite[Theorem 2]{F} we prove that  $Q_n$ and ${\bf q}_r=(q_{0,r},\ldots,q_{m-1,r})$, $(n,r)\in \mathbb{Z}_+^2$  satisfy the following  biorthogonality relation 
\begin{equation}\label{orthogonalityQ}
 \delta_{n,r}=\int Q_n(x) \sum_{j=0}^{m-1} q_{j,r}(x) \, d \sigma_j(x), \quad (n,r)\in \mathbb{Z}_+^2.
 \end{equation}
This relation in (\ref{orthogonalityQ}) allows us to call $\displaystyle Q_n$ and $\displaystyle {\bf q}_r$ type II and  type I  multiple orthogonal polynomials with respect to the distributions $(\sigma_0,\sigma_1, \ldots,\sigma_m)$, respectively. When $m=1$ both families coincide and they reduce to a sequence of standard orthogonal polynomials with respect to a single measure $\sigma_0$. Hence both type I and type II of multi-orthogonal polynomials are extensions of the regular orthogonal polynomials described in many publications such as \cite{AbSt, Chis, NiSo, stto}.  

During the first two decades of this century many works on type II multiple orthogonal polynomials on the starlike sets have been published (for instance \cite{LVA1, LVA2, AL1, AL2, AE}). They have two distinguishable starting points for the analysis in this topic. Some publications start by a higher order three terms of recurrence relations for polynomials to arrive to multi-orthogonality relations that these polynomials satisfy with respect to a system of measures supported on starlike sets. As examples of these papers we mention  \cite{AKS,AKV}. In the other staring point the authors take the multi-orthogonality of polynomials on a system of measures supported on the starlike sets and analyze the algebraic and asymptotic behavior of such polynomials. Apart of the mentioned before we find examples of these papers in \cite{BB, DL}. We now also analyze the type I multiple orthogonal polynomials o starlike sets. 

Given $(\alpha,\theta) \in \mathbb{R} \times [0, 2 \pi)$ set  $\mathcal{L} (\alpha,\theta)=(-\infty, \alpha]\exp i\theta \cup [\alpha, \infty)\exp i\theta$ the straight line with slope  $\tan \theta$ containing the point $\alpha$. Let us declare the following starlike sets
\begin{equation}\label{Sigmaset}
\begin{array}{l}
\displaystyle \Sigma_0 \qquad \mbox{as in (\ref{Sigma0})}, \\ \\
 \displaystyle \Sigma_e=\bigcup_{k=0}^m  \mathcal{L}\left(-\mu-\lambda, \frac{(2k+1) i \pi}{m+1}\right),\\ \\
\displaystyle \Sigma_o=\bigcup_{k=0}^m \left[-\mu-\lambda, +\infty\right) \exp \frac{2k i \pi}{m+1}.
\end{array}
\end{equation}
We denote
\[
\Omega_0=\mathbb{C}\setminus \Sigma, \quad \Omega_j=\mathbb{C}\setminus \{\Sigma_o\cup \Sigma_e\}, \quad  1 \leq j< m, 
\]
and
\[
\Omega_m=\left\{\begin{array}{l l}
\mathbb{C} \setminus \Sigma_e &  \mbox{if  } m \mbox{ even},\\ \\
 \mathbb{C} \setminus \Sigma_o &  \mbox{if } m \mbox{ odd}.
\end{array}\right.
\]
Let us choose the branch solutions $\omega_0,\, \omega_1,\ldots, \omega_m$ of the algebraic equation (\ref{omegaequationQ}) such that
\[
\omega_j\in \mathcal{H}\left(\Omega_j\right), \quad \omega_j(z)=\mathcal{O}\left(\frac{1}{z}\right) \quad \mbox{as} \quad z\to \infty, \quad j=1, \ldots, m,
\] 
and $\omega_0\in \mathcal{H}(\Omega_0)$ satisfying the asymptotic behavior stated in (\ref{asimptoticomega0}).

Let us introduce the algebraic equation 
\begin{equation}\label{omegaequationq}
\mu\, v^{m+1}-(z+\lambda+\mu)\, v+\lambda=0.
\end{equation}
Observe that the functions $v_j=1/\omega_j$, $j=0,1,\ldots,m$, are solutions of (\ref{omegaequationq}), where the functions $\omega_j$, $j=0,\ldots,m$, are the mentioned solutions of the algebraic equation (\ref{omegaequationQ}). The branch $v_m\in \Omega_m$ is such that 
\[
|v_m(z)|>|v_{m-1}|>\cdots >|v_0|, \quad z \in \Omega_m. 
\]

We are now ready to state our main results

\begin{theorem}\label{Ithaca} Let $\displaystyle \left\{Q_n\right\}_{n\in \mathbb{Z}_+}$ be the sequences of polynomials defined in (\ref{recurrenceQ}). The following equality holds uniformly for any compact set $K\subset \Omega_0$
\[
\lim_{n\to \infty} \frac{Q_n(z)}{\omega_0^n(z)}=\frac{1}{\lambda^m}\frac{1}{\omega_0^m(z)-\frac{m \mu}{\lambda \omega_0(z)}}\left(1+\left(-\frac{\mu}{\lambda}\right)^m\sum_{j=1}^{m} (-\mu\, \lambda \, \omega_0)^{j}(z)\right).
\]
\end{theorem}

The function $\displaystyle f_Q(z)=1+\left(-\frac{\mu}{\lambda}\right)^m\sum_{j=1}^{m} (-\mu\, \lambda \, \omega_0)^{j}(z)$ vanishes at a fixed number of points, and in general 
\[
\displaystyle \frac{1}{\lambda^m}\frac{1}{\omega_0^m(z)-\frac{m \mu}{\lambda \omega_0(z)}}\left(1+\frac{(-1)^{m+1}}{\lambda^m}\sum_{j=1}^{m} (-\lambda \, \omega)^{j}(z)\right) \in \mathcal{H}\left(\Omega_0\right)
\]
is holomorphic on $\Omega_0$.  Taking into account that $\Sigma_0$ is bounded, we deduce from Theorem \ref{Ithaca} that the zeros of the polynomials $Q_n$ are attracted by either the starlike set $\Sigma_0$ or by the zeros of the function $f_Q$. When we say that a set $A$ attracts another set $B$, we mean that $A$ contains all the cluster points of $B$. So the set of all zeros of the polynomials $Q_n$ is bounded.

\begin{theorem}\label{Ithaca2} Let $\displaystyle \left\{{\bf q}_r=(q_{0,r},\ldots,q_{m-1,r})\right\}_{r\in \mathbb{Z}_+}$ be the sequences of vector polynomials defined in (\ref{preccurencelatinq}). The following equality holds uniformly for any compact set $K\subset \Omega_m$
\[
\lim_{r\to \infty} \frac{{\bf q}_r}{v_m^r}=\frac{\displaystyle 1-\mu \left(\frac{v_m}{\lambda}\right)^m}{\displaystyle  \mu\, m v_m^{m+1}-1}\left(1, \frac{\lambda}{v_m}, \ldots, \left(\frac{\lambda}{v_m}\right)^{m-1}\right)
\]
\end{theorem}

Contrary to the previous case in Theorem \ref{Ithaca}, in Theorem \ref{Ithaca2} the starlike set $\Sigma_m$ is not bounded, hence the roots of polynomial components in ${\bf q}$ are attracted by either $\Sigma_m$ or infinity, and some of them could approach the zeros of the function $\displaystyle f_{\bf q}=1-\mu \left(\frac{v_m}{\lambda}\right)^m$. In this case the corresponding set of zeros is not necessarily bounded.

Theorem \ref{Ithaca} and Theorem \ref{Ithaca2} are proved in Section \ref{proof1} and Section \ref{proof2}, respectively. Previously, in Section \ref{auxiliar} we study some aspects about the distribution of zeros corresponding to certain polynomials envolved in the proofs of Theorem \ref{Ithaca} and Theorem \ref{Ithaca2}. This analysis justifies the reason why we choose the cuts in (\ref{Sigmaset}) for branches solutions of the algebraic equations in (\ref{omegaequationQ}) and (\ref{omegaequationq}).

\section{Types of generalizations of eigenvalue and eigenvector problems}\label{auxiliar}

In order to prove Theorem   \ref{Ithaca} and Theorem \ref{Ithaca2} in Section \ref{proof1} and Section \ref{proof2}, respectively, we use some results related to a variation of the eigenvalue and eigenvector problem. There are several previous publications (see for instance \cite{DL,DK}) on multi-orthogonality of polynomials that deal with certain generalizations of eigenvalues. We call partial eigenvalue and eigenvector problem to the variation that the present Section deals with.

Fix $n \in \mathbb{Z}_+$, and consider an $n \times n$ tridiagonal symmetric matrix $\displaystyle {\bf K}_n=[k_{i,j}]$, whose entries $k_{i,j}$, $(i,j)\in \{0,1,\ldots, n-1\}^2$ have the form:
\[
k_{i,j}=\left\{\begin{array}{l l ll}
b_j & \mbox{if} & j=i\pm1, ,\\  & & \\
a_j\not =0 & \mbox{if} & j=i, \\ & &  \\
0 & \mbox{otherwise}, & 
\end{array}
\right.
\]
where we requiere that $a_j\not =0$ either when $j=0$ or $j\not =(m+1) d_j$ with
\[
d_j=\left\lfloor \frac{j}{m+1}\right\rfloor, \quad j\in\{0,1,\ldots,n-1\}.
\]
The symbol $\displaystyle \left\lfloor \cdot \right\rfloor$ denotes the arithmetic floor function.

Set the $n \times n$ diagonal matrix function $\mathbb{I}_{n}(z)=i_{i,j}(z)$ as follows 
\[
i_{i,j}(z)=\left\{\begin{array}{r l l}
-z & \mbox{if}  & i=j \mbox{ and }  \frac{j}{m+1} \in \mathbb{N},  \\  &  &\\
0 & \mbox{otherwise}. &  
\end{array}
\right.
\]
Denote ${\bf S}_n(z)={\bf K}_n+\mathbb{I}_n(z)$. The value $x \in \mathbb{C}$ is said to be a partial eigenvalue of ${\bf K}_n$ if it is a root of the polynomial $\mathcal{P}(z)=\det {\bf S}_n(z)$.  By the structure of ${\bf K}_n$ we have that the polynomial $\mathcal{P}(z)$ has degree $\displaystyle d_n=\left\lfloor \frac{n}{m+1}\right\rfloor$. For each partial eigenvalue $x$ there are non zero column vectors ${\bf v}\in \mathbb{C}^n$ that satisfy ${\bf K}_n{\bf v}=\mathbb{I}_n(x) {\bf v}$. They are called partial eigenvectors.  

\begin{lemma}\label{manoinvisible} Fix $n=(m+1)\,d_n+k_n$ with $d_n\in \mathbb{Z}_+$ and $k_n\in \{0,1,\ldots,m-1\}$. The matrix ${\bf K}_n$ has $d_n$ different real partial eigenvalues. Equivalently, the polynomial $\mathcal{P}=\det {\bf S}_{n}$ has $d_n$ real simple zeros.
\end{lemma}

\begin{proof} First we prove that any partial eigenvalue is real. Recall that the matrix ${\bf K}_n$ is symmetric. Fix $x$ a partial eigenvalue of ${\bf K}_n$ and a corresponding partial eigenvector ${\bf v}$. There must be at least $j\in \{1,\ldots, d_n\}$ such that a component $v_{j\, (m+1)}\not =0$. Otherwise any complex number could be a partial eigenvalue of ${\bf K}_n$, and the polynomial $\mathcal{P}_n$ would be the zero constant function. 
 
Let us observe that ${\bf v}^{\dag}{\bf K}_n={\bf v}^{\dag}\mathbb{I}_n(\overline{x})$, where $\dag$ denotes the corresponding adjoint matrix, and $\overline{x}$ is the complex conjugate number of $x$. We have used that ${\bf K}_n$ and $\mathbb{I}$ are symmetric and real matrices. Finally we have that the modulus 
 \[
 \displaystyle ||{\bf K}_n {\bf v}||^2=||{\bf \bf v}||^2-(1- \overline{x})\sum_{j=1}^{d_n}||v_{(m+1)j}||^2=||{\bf \bf v}||^2-(1-x)\sum_{j=1}^{d_n}||v_{(m+1)j}||^2.
 \]
 Since not all the components $v_{j(m+1)}$, $j\in \{1,\ldots,d_n\}$ are zero, we conclude that $x=\overline{x}$, which proves that any partial eigenvalue $x$ is real.

 In order to prove that the partial eigenvalues are different we use the induction method. As an initial condition we observe that given $n \in \{m+1, \ldots, 2m+1\}$  $\deg \mathcal{P}_n=1$, hence the statement is trivially satisfied. We now assume that for an arbitrary $n\in \mathbb{N}$ the polynomial $\mathcal{P}_n$ has $d_n$ simple zeros, and then we prove that $\mathcal{P}_{n+1}$ vanishes $d_{n+1}$ different times. Let $x_{j,n}$, $j=1, 2, \ldots, d_n$ be the partial eigenvalues of ${\bf K}_n$ satisfying that 
 \[
 x_{1,n} <  x_{2,n}< \cdots <  x_{d_n,n}.
 \]
 By assumption, given $j \in \{1,\ldots, d_n\}$ the derivative
 \[
 \mathcal{P}_n^{\prime}(x_{j,n})=\lim_{z \to x_{j,n}} \frac{1}{z-x_{j,n}} \det {\bf S}_n(z) \not =0.
 \]
 Using the continuity of determinants we obtain that
 \[
 \lim_{z \to \lambda_{j,n}} (z-x_{j,n}) \det {\bf S}_n^{-1}(z)=\frac{1}{\mathcal{P}_n^{\prime}(x_{j,n})}.
 \]
 This implies that no column in the matrix function $(z-x_{j,n}){\bf S}_n^{-1}{\bf e}_k$, $k\in \{1,\ldots,n\}$ tends to the null vector. In this case ${\bf e}_k$ denotes the $k$th column in the $n \times n$ identity matrix. In particular, taking into account  that $\displaystyle \lim_{z\to \lambda_{j,n}}(z-x_{j,n}){\bf S}_n^{-1}(z)$ is a symmetric matrix, we have that 
 \begin{equation}\label{posittivedefinite}
 \lim_{z\to x_{j,n}}  {\bf e}_n^{\top}(z-x_{j,n}){\bf S}_n^{-1}(z) {\bf e}_n >0.
 \end{equation}

We now observe that the matrix function ${\bf S}_{n+1}(z)$ can be written by matrix blocks as follows
 \[
{\bf S}_{n+1}(z)=\left(\begin{array}{l r}
{\bf S}_n (z)& b_n{\bf e}_n \\ & \\
b_n {\bf e}_n^{\top} & c_n(z) 
\end{array}
\right) \quad \mbox{where} \quad c_n(z)=\left\{\begin{array}{l l}
a_n, & n\not=(m+1)d_n, \\ & \\
a_n-z, & n=(m+1)d_n.
\end{array}\right.
 \] 
 Hence 
 \[
 \mathcal{P}_{n+1}(z)= \det {\bf S}_n\,\, \det \left(c_n(z)- b_n^2{\bf e}_n^{\top}{\bf S}_n^{-1} {\bf e}_n\right)(z)=\mathcal{P}_n(z) \left(c_n(z)- b_n^2{\bf e}_n^{\top}{\bf S}_n^{-1} {\bf e}_n\right)(z).
  \]
Observe that $ \mathcal{P}_n$ and $ \mathcal{P}_{n+1}$ are  polynomials, then $b_n^2{\bf e}_n^{\top}{\bf S}_n^{-1} {\bf e}_n$ is a rational function with simple poles, which are the zeros of $\mathcal{P}_n$. Since $\deg \mathcal{P}_{n+1}=\deg \mathcal{P}_n+\deg c_n$ we have  the simple fraction decomposition
\[
\frac{\mathcal{P}_{n+1}}{\mathcal{P}_n}=c_n-b_n^2{\bf e}_n^{\top}{\bf S}_n^{-1} (z){\bf e}_n-=c_n+M-\sum_{j=1}^{d_n} \frac{\alpha_{j,n}}{z-x_{j,n}},
\]
 where 
 \[
M_n \in \mathbb{R} \quad \mbox{and} \quad  \alpha_{j,n}=\lim_{z\to x_{j,n}} (z-x_{j,n})\, b_n^2\, {\bf e}_n^{\top}{\bf S}_n^{-1} (z){\bf e}_n>0.
  \]
 In the last inequality we used the relation (\ref{posittivedefinite}). This implies that given two consecutive roots $x_{j,n}$ and $x_{j+1,n}$ of $\mathcal{P}_n$, the numerator $\mathcal{P}_{n+1}$ vanishes an odd number of times between them, hence the polynomial $\mathcal{P}_{n+1}$ has at least $d_n-1$ real zeros between $x_{1,n}$ and $x_{d_n,n}$. We observe that
 \[
 \lim_{z\to \lambda_{d_n,n}+} \frac{\mathcal{P}_{n+1}}{\mathcal{P}_n} =-\infty \quad \mbox{and}  \lim_{z\to \infty} \frac{\mathcal{P}_{n+1}}{\mathcal{P}_n}=\left\{\begin{array}{l c l}
 +\infty & \mbox{if} & d_{n+1}=d_n+1 \\ & & \\
 a_n+ M& \mbox{if} & d_{n+1}=d_n
 \end{array}
 \right. 
 \] 
 and
 \[
 \lim_{z\to \lambda_{1,n}-} \frac{\mathcal{P}_{n+1}}{\mathcal{P}_n} =\infty \quad \mbox{and}  \lim_{z\to \infty} \frac{\mathcal{P}_{n+1}}{\mathcal{P}_n}=\left\{\begin{array}{l c l}
 -\infty & \mbox{if} & d_{n+1}=d_n+1 \\ & & \\
 a_n+ M& \mbox{if} & d_{n+1}=d_n
 \end{array}
 \right. 
 \] 
Hence when $d_{n+1}=d_{n}+1$ the polynomial $\mathcal{P}_{n+1}$ must vanishes at least once before $x_{1,n}$ and once after $x_{d_n,n}$. This yields all the zeros of $\mathcal{P}_{n+1}$ are real. In the other case, depending of $a_n+M$'s sign the rational function vanishes either  before $x_{1,n}$ or after $x_{d_n,n}$. This completes the proof.
 \end{proof}
 
Let us consider a sequence of polynomials $\displaystyle \left\{h_n\right\}_{n\in \mathbb{Z}_+}$ whose elements are defined by the following initial conditions and recurrence relation:
\begin{equation}\label{initialhf}
h_0=h_1=\cdots=h_m=1,
\end{equation}
\begin{equation}\label{recurrencehf}
h_{n+1}(z)=\left\{\begin{array}{l c l}
h_n(z)-\mu \lambda^m h_{n-m}(z)  & \mbox{if} & \kappa < m, \\ & & \\
z \, h_n(z)-\mu \lambda^m h_{n-m}(z) & \mbox{if} & \kappa=m,
\end{array}
\right. \quad \gamma \geq 1.
\end{equation}
We write this recurrence relation in (\ref{recurrencehf}) with a matrix form 
\begin{equation}\label{recurrencehmatrix}
{\bf H}_n(z){\bf h}_n(z)=h_n(z){\bf e}_n, \quad n\in \mathbb{Z}_+,
\end{equation}
where  ${\bf h}_n$ is a vector column polynomial with $\displaystyle {\bf h}_n=\left(h_0\, h_1\, \cdots h_{n-1} \right)^{\top}$, ${\bf e}_n$ denotes the last column of the $n\times n$ identity matrix $\mathbb{I}_n$, and   ${\bf H}_n(z)=[h_{i,j}(z)] \in \mathbb{R}^{n\times n}$ is an $n\times n$ matrix  function whose entries $h_{i,j}\in \mathbb{R}$, $(i,j)\in \{0,1,2,\ldots,n-1\}^2$ satisfy:
\[
h_{i,j}(z)=\left\{\begin{array}{r c c}
-1 & \mbox{if} & j=i+1, \\ & & \\
z & \mbox{if}  & j=i=(m+1)d_j \\ & & \\
1 & \mbox{if}  & j=i\not=(m+1)d_j \\ & & \\
-\mu \lambda^m & \mbox{if} & i=j+m, \\ & & \\
0 & \mbox{otherwise}, &
\end{array}
\right. \quad \mbox{with} \quad d_j\in \mathbb{N}.
\]
After some  elementary matrix operations ${\bf H}_n$ is reduced to a tridiagonal matrix as ${\bf S}_n$, which means that the polynomials $\det {\bf H}(z)$ and $\det {\bf S}_n(z)$ are multiple. From the equality (\ref{recurrencehmatrix}) we deduce that if $h_n$ vanishes at a point $x$ then ${\bf H}_n(x) {\bf h}(x)=0$. This implies that the polynomials $\det {\bf H}(z)$, $\det {\bf S}_n(z)$, and $h_n$ are multiple. Therefore, taking into account Lemma \ref{manoinvisible}, we arrive at the following statement

\begin{proposition}\label{manoinvisibleH} Fix $n=m\,d_n+k_n$ with $d_n\in \mathbb{Z}_+$ and $k_n\in \{0,1,\ldots,m-1\}$. The polynomial $h_n$ has $d_n$ real simple roots.
\end{proposition}

Let us consider the polynomial elements of the sequence $\displaystyle \left\{p_r\right\}_{r\in \mathbb{Z}_+\cup \{-1\}}$ that satisfy the following initial conditions and recurrence relation:
\begin{equation}\label{initialh}
p_{-1}=0, \quad p_0=1, \quad p_1=\cdots=p_{m-1}=0,
\end{equation}
and when $r \geq m$
\begin{equation}\label{recurrenceh}
\mu \lambda^m\, p_{r}(z)=\left\{\begin{array}{l c l}
z \, p_{r-m}(z)+(-1)^m p_{r-m-1}(z) & \mbox{if} & \ell=0,\, k= 0,\\ & & \\
z \, p_{r-m}(z)+p_{r-m-1}(z) & \mbox{if} & \ell=0,\, k \not = 0,\\ & & \\
 p_{r-m}(z)+p_{r-m-1}(z)  & \mbox{if} & \ell \not =0, \, k\not =0,\\ & & \\
  p_{r-m}(z)+(-1)^mp_{r-m-1}(z)  & \mbox{if} & \ell \not =0, \, k=0.
\end{array}
\right. 
\end{equation}
This recurrence relation above has the matrix expression 
\begin{equation}\label{manoinvesibledual}
{\bf P}_r(z){\bf p}_r(z)=\mu \lambda^m p_r(z) {\bf e}_r, \qquad r\in \mathbb{Z}_+,
\end{equation}
where $\displaystyle {\bf p}_r=\left(p_0 \,\, p_1\,\,\cdots \,\, p_{r-1}\right)^{\top},$ $r\in \mathbb{Z}_+$ and the matrix function $\displaystyle {\bf P}_r=[p_{i,j}]\in \mathbb{R}^{r\times r}$ whose real entries satisfy
\[
p_{i,j}=\left\{\begin{array}{r c l}
-\mu \lambda^m & \mbox{if} & j=i+1, \\ & & \\ 
1 & \mbox{if} & i=j+m-1 \mbox{ and }  \frac{j}{m+1}  \not \in \mathbb{N},  \\ & & \\
z & \mbox{if} & i=j+m-1 \mbox{ and }  \frac{j}{m+1}  \in \mathbb{N} , \\ & & \\
(-1)^{m-1} & \mbox{if} & i=j+m \mbox{ and }  \frac{i}{m}  \in \mathbb{N},   \\ & & \\
1 & \mbox{if} & i=j+m  \mbox{ and }  \frac{i}{m}  \not \in \mathbb{N},  \\ & & \\
0 & \mbox{otherwise}. & 
\end{array}
\right. 
\]

Proceeding analogously as in the proof of Lemma \ref{manoinvisible}, we are able to extend its result to any tridiagonal symmetric $r \times r$ matrix function $\widetilde{\bf S}_r$ with $r=n+m-1$, such that it is obtained after tacking on $m-1$  columns and $m-1$ rows at the end of a matrix ${\bf S}_n$. This means that given $\displaystyle d_r=\left\lfloor \frac{r}{n}\right\rfloor$ and $k_r=r-nd_r$ we introduce $\displaystyle \tau_r=\left\lfloor \frac{d_r-k_r}{m+1}\right\rfloor$. The $\tau_r$ zeros of the polynomial $\det \widetilde{\bf S}(z)$ are real and simple. We also observe that we can obtain a matrix as $\widetilde{\bf S}_r$ making certain matrix elementary operations in ${\bf P}_r$, which implies that the polynomials $\det \widetilde{\bf S}_r(z)$ and $\det {\bf P}_r$ are multiple, and using the equality (\ref{manoinvesibledual}), we conclude that $p_r$ also vanishes at the same point as $\widetilde{\bf S}_r(z)$  has partial eigenvalues. 

\begin{proposition}\label{manoinvisibleP} Fix $r\in \mathbb{Z}_+$ with $\displaystyle d_r=\left\lfloor \frac{r}{m}\right\rfloor$ and $k_r=r-md_r$. Set also $\displaystyle \tau_r=\left\lfloor \frac{d_r-k_r}{m+1}\right\rfloor$. The polynomial $p_r$ has $\tau_r$ simple real roots. Agreeing  when $\tau_r=-1$, $p_r$  is the zero constant function.
\end{proposition}

\section{Proof of Theorem \ref{Ithaca}} \label{proof1}

Let us consider the change of variable $\displaystyle \zeta=z+\mu+\lambda$ to translate the starlike sets $\Sigma_0,$ $\Sigma_e$, and $\Sigma_o$ in (\ref{Sigmaset}) centered all at  $-\mu-\lambda$, to another collection of starlike sets $S_0$, $S_e$, $S_o$ whose centers are at $\zeta=0$ in the complex plane:
 \begin{equation}\label{Sset}
\begin{array}{l}
\displaystyle S_0=\bigcup_{k=0}^m \left[0, \frac{m+1}{m}\left(\frac{\mu\lambda^m}{m}\right)^{1/(m+1)} \right]\exp \frac{2 \pi i k}{m+1}, \\ \\
 \displaystyle S_e=\bigcup_{k=0}^m  \left(-\infty, +\infty\right) \exp \frac{(2k+1) i \pi}{m+1},\\ \\
\displaystyle S_o=\bigcup_{k=0}^m  \left(0, +\infty\right) \exp \frac{2k i \pi}{m+1}.
\end{array}
\end{equation}
Set the domains
\[
D_0=\mathbb{C}\setminus S_0, \quad D_j=\mathbb{C}\setminus \{S_o\cup S_e\}, \quad  1 \leq j< m, 
\]
and
\[
D_m=\left\{\begin{array}{l l}
\mathbb{C} \setminus S_e &  \mbox{if  } m \mbox{ even},\\ \\
 \mathbb{C} \setminus S_o &  \mbox{if } m \mbox{ odd}.
\end{array}\right.
\]

Set the functions $\xi_j(\zeta)=\lambda\omega_j(\zeta-\mu-\lambda)$, $j=0,1,\ldots,m$. Observe that they are $m+1$ different solution of the algebraic equation  (studied in \cite{AKS, AKV,F})
\begin{equation}\label{xiequationT}
 \xi^{m+1}-z\, \xi^m+\lambda^m \mu=0.
\end{equation}
For each $j\in \{0,1,2, \ldots,m\}$, $\displaystyle \xi_j \in \mathcal{H}\left(D_j\right)$, and 
\[
\xi_0(z)=z+\mathcal{O}(1) \quad \xi_j(z)=\mathcal{O}\left(\frac{1}{z^{1/m}}\right),  \quad \mbox{as} \quad z\to \infty, \quad j\in \{1,2,\ldots,m \}.
\] 
According to \cite[Proposition 1]{AKV}\label{solutionproperties} we have that
\begin{equation}\label{inequalitycomplexxi}
|\xi_0|\geq |\xi_1|\geq \cdots \geq |\xi_m|, \quad \mbox{on }\, \mathbb{C}.
\end{equation}
In $S_0$ the inequalities are strict. 

Then instead of considering the polynomials $Q_n$, $n\in \mathbb{Z}_+$ directly, we study an auxiliary family monic polynomials $L_n$ that we can obtain from $Q_n$ as follows
\[
L_n(\zeta)=\lambda^nQ_n(\zeta-\mu-\lambda), \qquad n\in \mathbb{Z}_+.
\]
Observe that they satisfy the following recurrence relation
\begin{equation}\label{recurrenceL}
L_{n+1}(\zeta)=\zeta L_n(\zeta)- \lambda^m\mu L_{n-m}(\zeta), \quad n \in \mathbb{Z}_+\setminus \{0,1,\ldots,m\},
\end{equation}
with initial conditions
\begin{equation}\label{recurrenceLi}
L_n(\zeta)=(\zeta-\mu)^n, \qquad n\in \{0,1,\ldots,m\}.
\end{equation}
 
The elements of the sequence $\displaystyle \left\{L_n\right\}_{n\in \mathbb{Z}_+}$ satisfy the same recurrence relation of the polynomials$\displaystyle \left\{T_n\right\}_{n\in \mathbb{Z}_+}$ 
\begin{equation}\label{recurrenceT}
T_{n+1}(\zeta)=\zeta T_n(\zeta)-\lambda^m \mu T_{n-m}(\zeta), \quad n \in \mathbb{Z}_+,
\end{equation}
with different initial conditions
\begin{equation}\label{recurrenceTi}
T_{-m}=T_{-m+1}=\cdots=T_{-1}=0 \quad \mbox{and} \quad T_0=1.
\end{equation}
 The sequence of polynomials $T_n$ is a particular case of a wide family of polynomials studied in several  publications, such as \cite{AKS,AKV, BB, DL}. In \cite{F} we include a proof that  given $n\in \mathbb{Z}_+$ with $n=(m+1) d_n+k_n$ where $d_n \in \mathbb{Z}_+$ and $k_n \in \{0,1,\ldots, m\}$, there is  a polynomial $h_n$ as the ones defined by the recurrence relation (\ref{recurrencehf}) and initial conditions (\ref{initialhf}) in Section \ref{auxiliar}, that satisfies that $T_n(z)=z^{\kappa} h_n(z^{m+1})$. By Proposition \ref{manoinvisibleH} we conclude that the zeros of $T_n$ lie in a starlike set that contains $S_0$.  Actually, if we consider the equality in (\ref{asymptoticT}) we are able to ensure that any cluster point corresponding to the set of these zeros belongs to  $S_0$.

For each $n\in \mathbb{Z}_+$ we write $T_n$ as follows
\[
T_n(z)=\left(a_0\xi_0^{n+m}+a_1\xi_1^{n+m}+\cdots+a_m \xi_m^{n+m}\right)(z), \quad z \in \bigcap_{j=0}^m D_j,
\]
whose coefficients satisfy the initial conditions
\[
a_0\xi_0^{n}+a_1\xi_1^{n}+\cdots+a_m \xi_m^{n}=\delta_{n,m}, \quad n \in \{0,1,2,\ldots,m\}.
\]
This is a system of linear equations whose unknown are the coefficients $a_j$, $j\in \left\{0,1,\ldots,m\right\}$. The corresponding matrix equation is ${\bf M}{\bf x}={\bf b}$ where
\[
{\bf M}=\left(\begin{array}{c c c c}
1 & 1 & \cdots & 1\\ & & & \\
\xi_0 & \xi_1 & \cdots & \xi_m \\ & & & \\
\vdots & \vdots & \ddots & \vdots \\ & & & \\
\xi_0^{m} &\xi_1^m & \cdots & \xi_m^m 
\end{array}
\right), \quad {\bf x}=\left(\begin{array}{c}
a_0 \\  \\
a_1 \\ \\
\vdots \\  \\
a_m
\end{array}
\right), \quad \mbox{and} \quad {\bf b}=\left(\begin{array}{c}
0 \\  \\
\vdots \\ \\
0 \\  \\
1
\end{array}
\right).
\] 
Using Cramer's rule we obtain that for each $\displaystyle j\in \left\{0, 1,\ldots, m\right\}$
\begin{equation}\label{coeficienta}
a_j=\frac{(-1)^{j}}{\displaystyle \prod_{k\not= j}(\xi_j-\xi_k)}.
\end{equation}
Observe that
\[
P(\zeta)= \prod_{k=1}^m(\zeta-\xi_k)=\zeta^{m+1}-z\zeta^m+\lambda^m \mu \quad P(\xi_j)=0.
\]
Taking its derivative we obtain that
\[
P^{\prime}(\zeta)=(m+1)\zeta^m-zm\zeta^{m-1}=\frac{m}{\zeta}\left(\zeta^{m+1}-z\zeta^m\right)+\zeta^m.
\]
Hence
\[
P^{\prime}(\xi_j)=\xi_j^m-\frac{m\lambda^m\mu}{\xi_j}=\frac{\xi_j^{m+1}-m\lambda^m\mu}{\xi_j}=\prod_{k\not= j}(\xi_j-\xi_k).
\]
We deduce 
\[
a_j=\frac{(-1)^j\xi_j}{\xi_j^{m+1}-m\lambda^m\mu}=\left(\frac{(-1)^j}{1-\frac{m \lambda^m\mu}{\xi_j^{m+1}}}\right) \frac{1}{\xi_j^{m}}, \qquad j\in \{0,1, \ldots, m\},
\]
which is the denominator in (\ref{coeficienta}). This implies that
\[
T_n(z)=\sum_{j=0}^m \frac{(-1)^j \xi_j^n(z)}{1-\frac{m \lambda^m \mu}{\xi_j^{m+1}(z)}}, \quad  z \in \bigcap_{j=0}^m D_j,
\]
and
\begin{equation}\label{formulaT}
\frac{T_n(z)}{\xi_0^n(z)}=\frac{1}{1-\frac{m \lambda^m \mu}{\xi_0^{m+1}(z)}}+\sum_{j=1}^m \frac{(-1)^j}{1-\frac{m \lambda^m \mu}{\xi_j^{m+1}(z)}}\frac{ \xi_j^n(z)}{ \xi_0^n(z)}, \quad  z \in \bigcap_{j=0}^m D_j.
\end{equation}
From the inequality (\ref{inequalitycomplexxi}), for each $j\in \{1,\ldots, m\}$, we can extend the equality
\begin{equation}\label{limavirreinal}
\lim_{n\to \infty} \left(\frac{ \xi_j}{ \xi_0}\right)^n=0, \quad \mbox{uniformly on any compact set} \quad K\subset D_0. 
\end{equation}
Combining (\ref{formulaT}) and (\ref{limavirreinal}) we arrive at
\begin{equation}\label{asymptoticT}
\lim_{n\to \infty} \frac{T_n}{\xi_0^n}=\frac{1}{1-\frac{m \lambda^m \mu}{\xi_0^{m+1}}} \quad \mbox{on} \quad  K\subset D_0.
\end{equation}

Since the function $\displaystyle \frac{1}{1-\frac{m \lambda^m \mu}{\xi_0^{m+1}}}$ is holomorphic on $D_0$ that never vanishes on it, the limit (\ref{asymptoticT}) implies that the starlike set $S_0$ attracts the roots of the polynomials $T_n$. 

Let us recall the polynomials $\displaystyle  \{L_n\}_{n\in \mathbb{Z}_+}$ defined by the recurrence relation (\ref{recurrenceL}) and initial conditions (\ref{recurrenceLi}).  We now write its $n$th element as linear combination of $T_k$, $k \leq n$. We observe previously that any arbitrary sequence of polynomials $\displaystyle \{K_n\}_{n\in \mathbb{Z}_+}$ which whose elements are linear combinations as
\[
K_n(x)=b_m T_n+b_{m-1}T_{n-1}+\cdots +b_0T_{n-m}, \quad b_{\ell} \in \mathbb{R}, \quad \ell\in \{0,1,\ldots,m\},
\]
satisfy the same recurrence relation as $\{L_n\}_{n\in \mathbb{Z}_+}$  holds in (\ref{recurrenceL}). If we now impose   the $L_n$'s initial conditions for the polynomials $K_n$ we have that $L_n=K_n$, $n\in \mathbb{Z}_+$. In order to obtain the proper coefficients $b_{\ell}\in \mathbb{R}$, $\ell=0,1,\ldots, n$, we solve the system of equations:
\[
\sum_{\ell=0}^m b_{\ell}T_{\ell+k}(x)=(x-\mu)^k, \quad k=0,\ldots,m.
\]
Since the coefficients $b_j$ are constant we fix $x=0$ to obtain the following linear system of equations
\[
b_0 T_{k}(0)+\cdots+b_mT_{k+m}(0)=(-\mu)^k, \quad k=0,1,\ldots,m.
\]
Taking into account that $T_n=x^n$, $n\in \{0,1,\ldots, m\}$, and using the recurrence relation (\ref{recurrenceT}), we obtain that
\[
b_0=1, \quad b_{m-k+1}(-\lambda^m \mu)=(-\mu)^k, \quad k=1,\ldots, m,
\]
or equivalently
\[
b_0=1, \quad b_j=\frac{(-\mu)^{m-j}}{\lambda^m}, \quad j=1,2,\ldots,m.
\]
Then
\[
L_n(x)=T_{n-m}+\left(-\frac{\mu}{\lambda}\right)^m\sum_{j=1}^{m} (-\mu)^{j} T_{n-m+j}.
 \]
Substituting (\ref{formulaT}), we have for any compact set $K\subset D_0$ with
\[
\lim_{n\to \infty} \frac{L_n(\zeta)}{\xi_0^n(\zeta)}=\lim_{n\to \infty}\left(\frac{T_{n-m}(\zeta)}{\xi_0^{n-m}(\zeta)} \frac{1}{\xi_0^m(\zeta)}+\left(-\frac{\mu}{\lambda}\right)^m\sum_{j=1}^{m} (-\mu)^{j} \frac{T_{n-m+j}(\zeta)}{\xi_0^{n-m+j}(\zeta)}\frac{1}{\xi_0^{m-j}(\zeta)}\right)
\]
\[
\lim_{n\to \infty} \frac{L_n(\zeta)}{\xi_0^n(\zeta)}=\frac{1}{\xi_0^m(\zeta)-\frac{m \lambda^m \mu}{\xi_0(\zeta)}}\left(1+\left(-\frac{\mu}{\lambda}\right)^m\sum_{j=1}^{m} (-\mu)^{j} \xi_0^{j}(\zeta)\right).
\]
We recall that $\zeta=z+\mu+\lambda$, and obtain that
\[
\lim_{n\to \infty} \frac{Q_n(z)}{\omega_0^n(z)}=\frac{1}{\lambda^m}\frac{1}{\omega_0^m(z)-\frac{m \mu}{\lambda \omega_0(z)}}\left(1+\left(-\frac{\mu}{\lambda}\right)^m\sum_{j=1}^{m} (-\mu\, \lambda \, \omega)^{j}(z)\right),
\]
which completes the proof. \hfill $\Box$

\section{Proof of Theorem \ref{Ithaca2}}\label{proof2}

Analogously to the previous Section \ref{Ithaca} we start by analyzing an auxiliary case. This is the sequence of vector polynomials
\begin{equation}\label{latin}
\displaystyle \mathcal{T}=\left\{{\bf t}_r=\left(t_{0,r}, t_{1,r}, \ldots, t_{m-1,r}\right)\right\}_{r\in \mathbb{Z}_+}
\end{equation}
generated by the initial conditions and  recurrence relation
\begin{equation}\label{preccurencelatin}
\begin{array}{l}
{\bf t}_{-1}=\left(0,0,\ldots, 0\right), \quad {\bf t}_{r}= {\bf e}_{r}^{\top}, \quad r \in \{0, 1, \ldots, m-1\}, \\ \\
\lambda^m \mu \,{\bf t}_r=x{\bf t}_{r-m}-{\bf t}_{r-m-1}, \quad r \in \mathbb{Z}_+\setminus \{0,1,\ldots, m-1\}.
\end{array}
\end{equation}
The sequence $\mathcal{T}$ is dual to the sequence $\displaystyle \left\{T_n\right\}_{n\in \mathbb{Z}_+}$ defined by the recurrence relation (\ref{recurrenceT}) and initial conditions (\ref{recurrenceTi}) with respect to the multi-orthogonality relations that the polynomials satisfy, and are expressed in \cite{AKS,AKV, BB, DL}. 

\begin{lemma}\label{eigenvalue}Consider the sequence of vector polynomials $\mathcal{T}$ as in (\ref{latin}) defined by the recurrence relation (\ref{preccurencelatin}).  For each $r \in \mathbb{Z}_+$ we have that
\begin{equation}\label{shift}
t_r=t_{0,r}=t_{1,r+1}=t_{2,r+2}=\cdots=t_{m-1,r+m-1}.
\end{equation}
\end{lemma}

\begin{proof}[Proof of Lemma \ref{eigenvalue}] Fix $j\in \{0,1,\ldots,m\}$. The recurrence relations (\ref{preccurencelatin}) can be written as
\begin{equation}\label{initialT0}
t_{j,r} \equiv 0, \quad r \in \{-1,0, \ldots, j-1, j+1, \ldots, m-1\} \quad \mbox{and} \quad t_{j,j}(z)=1,
\end{equation}
and
\[
\lambda^m \, \mu\, t_{j,r}=x\, t_{j,r-m}- t_{j,r-m-1}, \quad r \in \mathbb{Z}_+\setminus \{0,1,\ldots, m-1\}.
\]
We show that the expressions above are equivalent to the following
\begin{equation}\label{initialT}
t_{j,r} \equiv 0, \quad r \in \{-m+j-1, \ldots, j-1\} \quad \mbox{and} \quad t_{j,j}(z)=1,
\end{equation}
and
\begin{equation}\label{biquiniT}
\lambda^m\, \mu\,\, t_{j,r}=x\, t_{j,r-m}- t_{j,r-m-1}, \quad r \in \mathbb{Z}_+\setminus \{-m+j-1, \ldots, j-1\}.
\end{equation}
Given $\displaystyle j\in \left\{-1, 0,1,\ldots,m-1\right\}$, the conditions (\ref{initialT0}) and (\ref{initialT}) coincide in
\[
t_{j,r} = 0, \quad r \in \{-1,0, \ldots, j-1\} \quad \mbox{and} \quad t_{j,j}(z)=1.
\]
We declare that
\[
t_{j,r} = 0, \quad r \in \{-m+j-1, \ldots,-2\},
\]
which does not contradict (\ref{initialT0}). Taking into account (\ref{biquiniT}) we obtain that
\[
\lambda^m \, \mu\, t_{j,r}=x\, t_{j,r-m}- t_{j,r-m-1}=0, \quad r\in \{j+1, \ldots, m+j-1\}.
\]
This shows the equalities in (\ref{initialT}). Hence given $j\in \{0,\ldots, m-2\}$ the sequences $\displaystyle \left\{t_{j,r}\right\}_{r \in \mathbb{Z}_+}$ and $\displaystyle \left\{t_{j+1,r}\right\}_{r \in \mathbb{Z}_+}$ satisfy the same recurrence relation (\ref{biquiniT}) and according to (\ref{initialT}), they differ in the initial conditions with the shift 
\[
t_{j,r}=t_{j+1,r+1}=0, \quad r\in \{-m+j-1, \ldots, j-1, j+1,\ldots, m+j-1\}, 
\]
and  $t_{j,j}= t_{j+1,j+1}=1.$ This proves (\ref{shift}).
\end{proof}

From (\ref{shift}) we conclude that any sequence $\displaystyle \left\{t_{j,r}\right\}_{r\in \mathbb{Z}_+}$, $j\in \{0,1,\ldots, m-1\}$ contains all the information of the vector polynomials in $\displaystyle \left\{ {\bf t}_r\right\}_{r\in \mathbb{Z}_+}$. We focus on $\displaystyle \left\{t_{0,r}=t_r\right\}_{r\in \mathbb{Z}_+},$ whose elements, according to (\ref{preccurencelatin}), satisfy the recurrence relation
\begin{equation}\label{bellaciao}
\begin{array}{l}
t_{-1}=0, \quad t_{0}=1, \quad t_1=0, \ldots,  t_{m-1}=0,\\ \\
\lambda^m \, \mu \, t_{r}(z)=zt_{r-m}(z)-t_{r-m-1}, \quad r \in \mathbb{Z}_+.
\end{array}
\end{equation} 

Let us now consider the algebraic equation
\begin{equation}\label{lambdaequation}
\lambda^m \mu\, \phi^{m+1}(z)-z\phi(z)+1=0.
\end{equation}
Observe that
 \begin{equation}\label{omegalambda}
 \xi=1/\phi
 \end{equation}
 which is a solution of the equation (\ref{xiequationT}). The following statement can be easily shown considering (\ref{omegalambda}).  There is a multivalued solution $\phi=(\phi_0,\phi_1, \ldots, \phi_m)$ of the equation (\ref{lambdaequation}). Note this equation is an algebraic function of order $m+1$ with global branches of $\phi(z)$ in (\ref{lambdaequation}) satisfying for each $\displaystyle j\in \left\{0,1,\ldots,m\right\}$, that $\phi_j \in \mathcal{H}(D_j)$
and 
\[
|\phi_0| < |\phi_1| < \cdots < |\phi_m|, \quad \mbox{on }\, D=\bigcup_{j=0}^m D_j.
\]
They also behave asymptotically 
\[
\begin{array}{l}
\displaystyle \phi_0(z)=\frac{1}{z} +\mathcal{O}\left(\frac{1}{z^2}\right) \quad \mbox{as} \, z \to \infty, \\ \\
\displaystyle \lim_{z\in D,\, z\to \infty} \frac{\lambda^m \mu \, \phi_j^m(z)}{z} =1, \quad  j=1,\ldots,m. 
\end{array}
\]

\begin{proposition}\label{algebraic} Set $r\in \mathbb{Z}_+$ with $r=d_r m+k_r $, $d \in \mathbb{Z}_+$ and $k_r \in \{0,1,\ldots, m-1\}$. Take $d_r-k_r=(m+1) \tau_r+\ell_r$  with $\tau_r \in \mathbb{Z}_+\cup \{-1\}$ and $\ell_r \in \{0,1,\ldots, m\}$.  There exists a polynomial $p_r$ with a nonnegative leading coefficient and degree $\tau_r$ such that 
\begin{equation}\label{treduced}
t_r(z)=\left\{\begin{array}{l c l}
(-1)^k_r z^{\ell_r} p_r(z^{m+1}) & \mbox{if} & \ell_r \in \{0,1,\ldots,m \}, \\ & & \\
0 & \mbox{if} & \ell_r \in \{-m+1, -m+2, \ldots, -1 \}.
\end{array}
\right.
\end{equation}
Agreeing $\deg p_r=-1$ means $p_r \equiv 0$. The polynomial elements of the sequence $\displaystyle \left\{p_r\right\}_{r\in \mathbb{Z}_+}$ can be constructed by the combination of the initial conditions in (\ref{initialh}) and the recurrence relation in (\ref{recurrenceh}). The polynomial $t_r$ has all its zeros in a starlike set that contains $S_m=\mathbb{C}\setminus D_m$.
\end{proposition}

\begin{proof}[Proof of Lemma \ref{algebraic}] According to the initial conditions in the recurrence relations (\ref{preccurencelatin}) we have that ${\bf t}_0=\left(1,0,\ldots,0\right)$. Taking into account (\ref{shift}) we obtain that $t_{0}=1=p_0$ and $t_j=0$, $j=1,\ldots,m-1$, which proves (\ref{initialh}). When $r=0$, $d_r=k_r=0$, hence $\tau_r=0$ and $\ell_r=0$, which is coherent with the fact $\deg t_0=0$. We observe that if $r \in \{1, \ldots, m-1\}$, then $d_r=0$, $r=k_r$ and $d_r-k_r=-k_r=-(m+1)+\ell_r$, hence $\tau_r=-1$, which agrees with $p_r \equiv 0$, $r\in \{1, \ldots, m-1\}$. 

In oder to prove (\ref{recurrenceh}) we use the induction method. Fix $r=d_rm+k_r \geq m$ which implies that $r-m=(d_r-1)m+k_r$, $k_r\in \{0,1,\ldots, m-1\}$ and 
\[
r-m-1=\left\{\begin{array}{l c l}
(d_r-2)m+m-1& \mbox{if} &k_r=0,\\ & & \\
(d_r-1)m+k_r-1 & \mbox{if} & k_r\in \{1,\ldots, m-1\}.
\end{array}
\right.
\]
Set $d_r-k_r=(m+1)\tau_r+\ell_r$.  We obtain that 
\[
(d-1)-k_r=\left\{\begin{array}{l c l}
(m+1)(\tau_r-1)+m & \mbox{if} &  \ell_r=0, \\ & & \\
(m+1)\tau_r+\ell_r-1 & \mbox{if} &  \ell_r \in \{1,\ldots,m\}.
\end{array}
\right.
\]
If $k_r\not =0$ we have that $(d_r-1)-(k_r-1)=d_r-k_r=(m+1)\tau_r+\ell_r$, and when $k_r=0$, we find that $d_r-2-m+1=d_r-(m+1)=(m+1)(\tau_r-1)+\ell_r,$ $\ell_r \in \{0,1,\ldots, m\}$.

Set the induction hypothesis
\[
t_{r-m}(z)=\left\{\begin{array}{l l l}
(-1)^{k_r} z^{m} p_{r-m}(z^{m+1}), & \deg p_{r-m}=\tau_r-1, & \ell_r=0, \\ & & \\ 
(-1)^{k_r} z^{\ell_r-1} p_{r-m}(z^{m+1}), & \deg p_{r-m}=\tau_r, & \ell_r \in \{1,\ldots, m\}.
\end{array}
\right.
\]
and 
\[
t_{r-m-1}(z)=\left\{\begin{array}{l c l}
(-1)^{m-1} z^{\ell_r} p_{r-m-1}(z^{m+1})& \mbox{if} & k_r=0, \\ & & \\
 (-1)^{k_r-1} z^{\ell_r} p_{r-m-1}(z^{m+1})& \mbox{if} & k_r\not=0,
 \end{array}
 \right. \quad \ell_r \in \{0,1,\ldots,m\}.
\]
\[
 \deg p_{r-m-1}=\left\{\begin{array}{l l l}
\tau_r-1 & \mbox{if} &  k_r=0, \\ & & \\ 
\tau_r & \mbox{if} &  k_r\not=0. 
\end{array}
\right.
\]
We divide the proof in cases considering different combinations of $k_r$ and $\ell_r$.
\begin{enumerate}
\item[]
\item {\bf Case 1 $k_r=0$ and $\ell_r=0$:} According to (\ref{bellaciao}) we obtain the following chain of equalities
\[
\lambda^m \mu\, t_{r}(z)=zt_{r-m}(z)-t_{r-m-1}=
\]
\[
=z z^m p_{r-m}(z^{m+1})+(-1)^{m}p_{r-m-1}(z^{m+1})
\]
\[
=z^{m+1} p_{r-m}(z^{m+1})+(-1)^m p_{r-m-1}(z^{m+1}).
\]
Recall that $t_r (z)=p_r(z^{m+1})$. Then 
\[
\lambda^m\mu \, p_r(z)=zp_{r-m}(z)+(-1)^mp_{r-m-1}(z).
\]
Taking into account  $\deg p_{r-m}=\deg p_{r-m-1}=\tau_r-1$  we have that 
\[
\deg p_r=1+\deg p_{r-m}=1+\tau_r-1=\tau_r. 
\]
Hence the leading coefficient of $p_r$ is positive.
\item[]
\item {\bf Case 2 $k_r>0$ and $\ell_r=0$:} Consider (\ref{bellaciao}) to obtain that
\[
\lambda^m \mu\, t_{r}(z)=zt_{r-m}(z)-t_{r-m-1}=
\]
\[
=(-1)^{k_r} z z^m p_{r-m}(z^{m+1})+(-1)^kp_{r-m-1}(z^{m+1})
\]
\[
=(-1)^{k_r}(z^{m+1} p_{r-m}(z^{m+1})+p_{r-m-1}(z^{m+1})).
\]
We observe that we can write $t_r (z)=(-1)^k p_r(z^{m+1})$ which satisfies that 
\[
\lambda^m\mu \, p_r(z)=zp_{r-m}(z)+p_{r-m-1}(z).
\]
Observe that $\deg p_{r-m}=\tau_r-1$ and $\deg p_{r-m-1}=\tau_r$.  Taking into account that the leading coefficients of $p_{r-m}$ and $p_{r-m-1}$ are both positive we have that $\deg p_r=1+\deg p_{r-m}=1+\tau_r-1=\tau_r=\deg p_{r-m-1}$. This implies that $p_r$ also has a positive leading coefficient.
\item[]
\item {\bf Case 3 $k_r=0$ and $\ell_r>0$:} Let us consider (\ref{bellaciao})  
\[
\lambda^m \mu\, t_{r}(z)=zt_{r-m}(z)-t_{r-m-1}=
\]
\[
=z z^{\ell_r-1} h_{r-m}(z^{m+1})+(-1)^{m}z^{\ell}h_{r-m-1}(z^{m+1})
\]
\[
=z^{\ell_r} \left[ p_{r-m}(z^{m+1})+(-1)^{m}p_{r-m-1}(z^{m+1})\right]
\]
Then we have that $t_r=z^{\ell_r} p_{r}(z^{m+1})$, with
\[
 \lambda^m\mu\, p_{r}(z)=p_{r-m}(z)+(-1)^{m}h_{r-m-1}(z)
\]
Since $\deg p_{r-m}=\tau_r$ and $\deg p_{r-m-1}=\tau_r-1$, the polynomial $p_r$ has degree $\tau_r$ with positive coefficient. 
\item[]
\item {\bf Case 4 $k_r>0$ and $\ell_r > 0$:} In this last case we also star by considering the recurrence relation (\ref{bellaciao}), hence
\[
\lambda^m \mu\, t_{r}(z)=zt_{r-m}(z)-t_{r-m-1}
\]
\[
=(-1)^{k_r}\,z \,z^{\ell-1}p_{r-m}(z^{m+1})+(-1)^{k_r} z^{\ell} p_{r-m-1}(z^{m+1})
\]
\[
=(-1)^{k_r} z^{\ell_r} \left[p_{r-m}(z^{m+1})+p_{r-m-1}(z^{m+1})\right],
\]
hence $t_r=(-1)^{k_r}z^{\ell}p_r(z^{m+1})$ with
\[
\lambda^m \mu\, p_{r}(z)= p_{r-m}(z)+p_{r-m-1}(z).
\]
Since $\deg p_{r-m}=\deg p_{r-m-1}=\tau_r$, then $\deg p_r=\tau_r$  with a positive leading coefficient.
\item[]  
\end{enumerate}
This completes the proof of Proposition \ref{algebraic}.
\end{proof}

Assuming (\ref{treduced}) the statement that any zero of $t_r$ lies in a starlike set that contains $S_m=\mathbb{C}\setminus D_m$ is consequence of Proposition \ref{manoinvisibleP} in Section \ref{auxiliar}. Furthermore, a consequence of Proposition \ref{strongbehavior} below is that $S_m$ contains all the cluster points that the set of zeros corresponding to the polynomials $t_r$. 

\begin{proposition}\label{strongbehavior} Consider the sequence $\displaystyle \left\{t_r\right\}_{r\in \mathbb{Z}_+}$ whose elements are the polynomials in (\ref{shift}). Then
\begin{equation}\label{acumulacion}
\lim_{r\to \infty} \frac{t_r}{\phi_m^r} =  \frac{\displaystyle 1}{\displaystyle \lambda^m \mu m \phi_m^{m+1}-1}\in \mathcal{H}\left(D_m\right),
\end{equation}
uniformly in any compact subset of $D_m$.
\end{proposition}

Since  the functions $\phi_m$ is a solution of the algebraic equation (\ref{lambdaequation}), it  does not vanish in $D_m$, hence the equality (\ref{acumulacion}) in Theorem \ref{strongbehavior} implies the set $S_m$ attracts the roots of the polynomials $t_r$, as $r \to \infty$ (considering that $S_m$ contains the infinity). A precedent of this result can be found in \cite{LS}. In the context of Nikishin systems of measures supported on the real line (see \cite{nik}), the authors proved that all the zeros of type I polynomials are lying in a set which is analogous to our $S_m$. They use the structure of Nikishin systems to prove this statement about the zeros distribution.

\begin{proof}[Proof of Proposition \ref{strongbehavior}]
Let us observe that the following function satisfies the recurrence relation (\ref{bellaciao}) and consider the distribution of zeros  described in Proposition \ref{algebraic}. This justifies the equality 
\begin{equation}\label{loureiro}
t_r= b_0  \phi_0^{r+m}+b_1 \, \phi_1^{r+m}+\cdots+b_m \phi_m^{r+m}, \quad r\in \mathbb{Z}_+.
\end{equation}
The coefficients $b_0(z),\ldots,b_m(z)$ are defined from the initial conditions
\[
t_r=0, \quad r=-m,-m+1, \ldots, -1, \quad \mbox{and} \quad t_0=1,
\]
which implies that
\[
\left\{\begin{array}{l l}
b_0+b_1+b_2+\cdots+b_m& =0,\\ & \\
\displaystyle b_0 \phi_0+b_1  \phi_1+b_2 \phi_2+\cdots+b_m \phi_m& =0,\\& \\
\cdots \cdots \cdots \cdots \cdots \cdots & \\ & \\
\displaystyle b_0 \phi_0^{m-1}+b_1 \phi_1^{m-1} +b_2 \phi_2^{m-1}+\cdots+b_m \phi_m^{m-1}& =1.
\end{array}\right.
\] 
Set ${\bf F}\, {\bf b}={\bf e}_m$ the matrix expression of the linear system of equations above, where
\[
{\bf F}=\left(\begin{array}{c c c c}
1 &  1 & \cdots &  1  \\ 
\phi_0 &  \phi_1 & \cdots &  \phi_m  \\
\vdots & \vdots & \ddots & \vdots \\
\phi_0^{m-1} &  \phi_1^{m-1} & \cdots &  \phi_m^{m-1}
\end{array}
\right) \quad 
{\bf b}=\left(\begin{array}{ c}
b_{0}  \\
b_{1}  \\
 \vdots \\
 b_{m} 
\end{array}
\right), \quad \mbox{and} \quad {\bf e}_m=\left(\begin{array}{ c}
0  \\
\vdots \\
 0 \\
1
\end{array}
\right).
\]
From Cramer's rule we find
\begin{equation}\label{Aj}
b_j=\frac{1}{\displaystyle \phi_j \prod_{k=0,k\not=j}^m (\phi_j-\phi_k)}, \qquad j=0,\ldots,m.
\end{equation}
Let us denote the polynomial in $\zeta \in \mathbb{C}$:
\[
R_z(\zeta)=\lambda^m \mu \prod_{k=0}^m (\zeta-\phi_k)=\lambda^m \mu \zeta^{m+1}-z\,\zeta+1,
\]
hence for each $\displaystyle j \in \left\{0,1,\ldots,m\right\}$
\[
R_z^{\prime}\left(\phi_j\right)= \lambda^m \mu \prod_{k=0,k\not=j}^m \left(\phi_j-\phi_k\right)=\lambda^m \mu\, (m+1)\, \phi_j^{m}-z.
\]
Since 
\[
R_z(\phi_j)=\lambda^m \mu \prod_{k=0}^m (\phi_j-\phi_k)=\lambda^m \mu\phi_j^{m+1}-z\,\phi_j+1=0,
\]
we have that
\begin{equation}\label{derivativenonzero}
\phi_j R_z^{\prime}\left(\phi_j\right)=\lambda^m \mu\, (m+1)\, \phi_j^{m+1}-z\phi_j=\lambda^m \mu m \phi_j^{m+1}-1.
\end{equation}
The $\phi_j R_z^{\prime}\left(\phi_j\right)$ never vanishes in $D_j$, this proves that
\[
\frac{1}{\displaystyle \lambda^m \mu m\phi_m^{m+1}-1}\in \mathcal{H}\left(\mathbb{C}\setminus S_m\right).
\] 
We substitute the expression (\ref{derivativenonzero}) in (\ref{Aj})
\[
b_j=\frac{1}{\displaystyle \lambda^m \mu m \phi_j^{m+1}-1}, \qquad j=0,\ldots,m.
\]
Considering (\ref{loureiro}) we obtain that
\[
t_r=\sum_{j=0}^{m}  \frac{\displaystyle \phi_j^{r+m}}{\displaystyle \lambda^m \mu m \phi_j^{m+1}-1}.
\]
Then
\[
\frac{t_r}{\phi_m^r} = \sum_{j=0}^{m}  \frac{\displaystyle 1}{\displaystyle \lambda^m \mu m \phi_j^{m+1}-1} \left(\frac{\phi_j}{\phi_m}\right)^{r-1},
\]
hence
\[
\lim_{r\to \infty} \frac{t_r}{\phi_m^r} = \frac{\displaystyle 1}{\displaystyle \lambda^m \mu m \phi_m^{m+1}-1} ,
\]
uniformly in any compact subset of $\mathbb{C}\setminus S_m$. The proof of  Proposition \ref{strongbehavior} is then completed.
\end{proof}

Combining Lemma \ref{eigenvalue} and Proposition \ref{strongbehavior} we obtain the following result 

\begin{corollary}\label{vectortas} Consider the sequence of vector polynomials 
\[
\displaystyle \displaystyle \mathcal{T}=\left\{{\bf t}_r=\left(t_{0,r}, t_{1,r}, \ldots, t_{m-1,r}\right)\right\}_{r\in \mathbb{Z}_+}
\] 
defined in (\ref{latin}) by  the  recurrence relation in (\ref{preccurencelatin}). Then 
\[
\lim_{r\to \infty} \frac{1}{\phi_m^r}{\bf t}_r = \frac{\displaystyle 1}{\displaystyle \lambda^m \mu m \phi_m^{m+1}-1}\left(1, \frac{1}{\phi_m}, \ldots, \frac{1}{\phi_m^{m-1}}\right) ,
\]
uniformly in any compact subset of $\mathbb{C}\setminus S_m$.
\end{corollary}

Let us introduce the sequence of of vector polynomials $\displaystyle \left\{{\bf u}_r\right\}$ with elements
\begin{equation}\label{JeanPierre}
{\bf u}_r(z)=\left\{\begin{array}{l c l}
{\bf q}_r & \mbox{if} & r \in \{-1,0,1,\ldots, m-1\},\\ & & \\
\displaystyle \frac{1}{\lambda^{r}} {\bf q}_r(z-\lambda-\mu) & \mbox{if} &  r\geq m,
\end{array}
 \right.
\end{equation}
whose elements satisfy the recurrence relation
\begin{equation}\label{preccurencelatinu}
\begin{array}{l}
{\bf u}_{-1}=\left(0,0,\ldots, 0\right), \quad {\bf u}_{r}= {\bf e}_{r}^{\top}, \quad r \in \{0, 1, \ldots, m-1\}, \\ \\
\lambda^m \mu \,{\bf u}_r(z)=(z-\mu){\bf u}_{r-m}(z)-{\bf u}_{r-m-1}(z),  \quad r \in \{m, \ldots, 2m-1\}, \\ \\
\lambda^m \mu \,{\bf u}_r(z)=z{\bf u}_{r-m}(z)-{\bf u}_{r-m-1}(z), \quad r \geq 2m.
\end{array}
\end{equation}

Since the recurrence relations of the vector polynomials $\{{\bf t}_r\}$ and $\{{\bf u}_r\}$ have the same structure when $r \geq 2m$, we write ${\bf u}_r$ in terms a linear form with elements of $\{{\bf t}_r\}$:
\[
{\bf u}_r=\beta_{2m-1} {\bf t}_r+\beta_{2m-2} {\bf t}_{r-1}+\cdots +\beta_{0} {\bf t}_{r-2m+1}, \quad \beta_{\ell}\in \mathbb{R}, \quad \ell\in \{0,1, \ldots,2m-1\}.
\]
We find the coefficients $\beta_{\ell}\in \mathbb{R},$ $\ell\in \{0,1, \ldots,2m-1\}$ using the initial conditions 
\[
\begin{array}{l l}
\displaystyle \sum_{\ell\geq 2m-r-1}^{2m-1} \beta_{\ell} {\bf t}_{\ell+r-2m+1}={\bf e}_r^{\top}, & r=0, \ldots,m-1,  \\ & \\
\displaystyle \lambda^m\mu \sum_{\ell=m-r-1}^{2m-1} \beta_{\ell} {\bf t}_{\ell+r-m+1}=(z-\mu){\bf e}_r^{\top}-{\bf e}_{r-1}^{\top}, & r=0, \ldots,m-1.  
\end{array}
\]
According to the initial conditions in (\ref{preccurencelatin}), the first group of equations becomes that
\[
 \sum_{\ell= 2m-r-1}^{2m-1} \beta_{\ell} {\bf e}^{\top}_{\ell+r-2m+1}={\bf e}_r^{\top}, \quad  r=0, \ldots,m-1,
\]
which implies that $\beta_{2m-1}=1$ and $\beta_{\ell}=0$, $\ell=m, \ldots, 2m-2$. We now consider the second group of equalities, that becomes
\[
\lambda^m\mu \left({\bf t}_{r+m}+\sum_{\ell=m-r-1}^{m-1} \beta_{\ell} {\bf t}_{\ell+r-m+1}\right)=(z-\mu){\bf e}_r^{\top}-{\bf e}_{r-1}^{\top}, \quad  r=0, \ldots,m-1.
\]
Using again the initial conditions  (\ref{preccurencelatin}) we arrive at
\[
z{\bf e}^{\top}_r-{\bf e}^{\top}_{r-1}+\sum_{\ell=m-r-1}^{m-1} \beta_{\ell} {\bf e}_{\ell+r-m+1}^{\top}=(z-\mu){\bf e}_r^{\top}-{\bf e}_{r-1}^{\top}, \quad r=0, \ldots,m-1.
\]
Hence $\beta_{m-1}=-\mu$ and $\beta_{\ell}=0$, $\ell=0,1,\ldots, m-2$, and 
\[
{\bf u}_r={\bf t}_r-\mu{\bf t}_{r-m}, \qquad r \in \mathbb{Z}_+,
\]
or equivalently
\[
\frac{{\bf u}_r}{\phi_m^r}=\frac{{\bf t}_r}{\phi_m^r}-\mu\phi_m^m\frac{{\bf t}_{r-m}}{\phi_m^{r-m}}, \qquad r \in \mathbb{Z}_+.
\]
Combining this identity and Corollary \ref{vectortas} we obtain the following result

\begin{corollary}\label{vectortasu} Consider the sequence of vector polynomials 
\[
\displaystyle \left\{{\bf u}_r=\left(u_{0,r}, u_{1,r}, \ldots, u_{m-1,r}\right)\right\}_{r\in \mathbb{Z}_+}
\] 
whose elements satisfy the relations in  (\ref{preccurencelatinu}). Then 
\[
\lim_{r\to \infty} \frac{1}{\phi_m^r}{\bf u}_r = \frac{\displaystyle 1-\mu \phi_m^m}{\displaystyle \lambda^m \mu m \phi_m^{m+1}-1}\left(1, \frac{1}{\phi_m}, \ldots, \frac{1}{\phi_m^{m-1}}\right) ,
\]
uniformly in any compact subset of $\overline{\mathbb{C}}\setminus S_m$.
\end{corollary}
The proof of Theorem \ref{Ithaca2} is completed recovering the change of variable in (\ref{JeanPierre}) and taking into account that  the solutions $v_0, v_1, \ldots, v_m$ of the algebraic equations (\ref{omegaequationq}) are related with the solutions $\phi_0, \phi_1, \ldots, \phi_m$ of (\ref{lambdaequation}) by the equalities
\[
v_j(z)=\lambda\phi_j(z+\lambda+\mu), \quad j=0,1,\ldots,m.
\] 

\bibliographystyle{amsplain}

\end{document}